\documentclass[10pt]{amsart}
\usepackage{amssymb, amsthm, amsmath}
\usepackage{xy}
\usepackage{epsfig}
\usepackage{psfrag}

\psfrag{lambda}{$\lambda$}
\theoremstyle{plain}
\newtheorem{prop}{Proposition}
\newtheorem{thm}{Theorem}
\newtheorem{cor}{Corollary}
\newtheorem{lemma}{Lemma}
\newtheorem{conj}{Conjecture}
\theoremstyle{remark}
\theoremstyle{definition}

\newtheorem{defn}{Definition}
\newtheorem{exa}{Example}

\newcommand{\CC}{{\mathcal M}}
\newcommand{\II}{{\mathcal I}}
\newcommand{\bP}{{\mathbb P}}

\newcommand{\Z}{{\mathbb Z}}

\newcommand{\C}{{\mathbb C}}
\newcommand{\D}{{\mathcal D}}
\newcommand{\cO}{{\mathcal O}}
\newcommand{\E}{{\mathcal E}}

\newcommand{\zZ}{{\mathcal Z}}
\newcommand{\W}{{\mathcal W}}
\newcommand{\s}{\sigma}
\newcommand{\la}{\lambda}
\newcommand{\m}{\mu}
\newcommand{\n}{\nu}
\newcommand{\cc}{\negmedspace:\negmedspace}
\newcommand{\bull}{{\scriptscriptstyle \bullet}}

\DeclareMathOperator{\Span}{Span}
\DeclareMathOperator{\Ker}{Ker}

\newcommand{\dis}{\displaystyle}

\newcommand{\ssm}{\smallsetminus}
\newcommand{\gequ}{\geqslant}
\newcommand{\lequ}{\leqslant}

\newcommand{\ra}{\rightarrow}
\newcommand{\skipline}{\vspace{\baselineskip}}

\newcommand{\ov}{\overline}
\newcommand{\noin}{\noindent}
\newcommand{\wt}{\widetilde}
\newcommand{\wh}{\widehat}

\newcommand{\pzp}[1]{\includegraphics[scale=.9]{#1.ps}}
\newcommand{\pzs}{\hspace{4mm}}

\begin{document}

\title[Gromov-Witten invariants and quantum cohomology]
{Gromov-Witten invariants and quantum cohomology of Grassmannians}
\author{Harry Tamvakis}
\date{June 27, 2003}
\address{Department of Mathematics, Brandeis University - MS 050,
P. O. Box 9110, Waltham, MA
02454-9110, USA}
\email{harryt@brandeis.edu}

\begin{abstract} 
This is the written version of my five lectures at the Banach Center 
mini-school on `Schubert Varieties', in Warsaw, Poland, May 18--22,
2003.
\end{abstract}

\maketitle

\section{Lecture One}
\label{l1}

\noindent
The aim of these lectures is to show that three-point genus zero Gromov-Witten
invariants on Grassmannians are equal (or related) to classical triple
intersection numbers on homogeneous spaces of the same Lie type, and to use
this to understand the multiplicative structure of their (small) quantum
cohomology rings.  This theme will be explained in more detail as the lectures
progress. Much of this research is part of a project with Anders S.\ Buch
and Andrew Kresch, presented in the papers \cite{Buch}, \cite{KTlg},
\cite{KTorth}, and \cite{BKT1}. I will attempt to give the original references
for each result as we discuss the theory.

\subsection{The classical theory}
We begin by reviewing the classical story for the type $A$ Grassmannian.
Let $E=\C^N$ and $X=G(m,E)=G(m,N)$ be the Grassmannian of $m$-dimensional
complex linear subspaces of $E$. One knows that $X$ is a smooth projective 
algebraic variety of complex dimension $mn$, where $n=N-m$. 

\begin{figure}[htpb]
\epsfxsize 35mm
\center{\mbox{\epsfbox{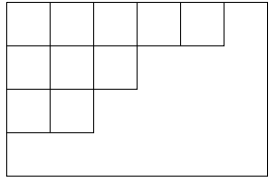}}}
\caption{The rectangle $(6^4)$ containing $\la=(5,3,2)$}
\label{rectangle}
\end{figure}

The space $X$ is stratified by Schubert cells; the closures of these
cells are the {\em Schubert varieties} $X_{\la}(F_\bull)$,
where $\la$ is a partition and 
\[
F_{\bull} \ :\ 0=F_0\subset F_1\subset \cdots \subset F_N=E
\]
is a complete flag of subspaces of $E$, with $\dim F_i=i$ for each $i$.
The partition $\la =(\la_1\gequ \la_2 \gequ \cdots \gequ \la_m \gequ 0)$
is a decreasing sequence of nonnegative integers such that $\la_1\lequ n$. 
This means that the Young diagram of $\la$ fits inside an $m\times n$
rectangle, which is the diagram of $(n^m)$. We denote this containment
relation of diagrams by $\la\subset (n^m)$. The diagram 
shown in Figure \ref{rectangle} corresponds to a Schubert
variety in $G(4,10)$.

The precise definition of $X_{\la}(F_\bull)$ is 
\begin{equation}
\label{defequ}
 X_\lambda(F_\bull) = \{ \,V \in X \mid \dim(V \cap
   F_{n+i-\lambda_i}) \gequ i, ~\forall \ 1 \lequ i \lequ m \,\} \,.
\end{equation}
Each $X_{\la}(F_\bull)$ is a closed subvariety of $X$ of codimension equal to
the weight $|\lambda| = \sum \lambda_i$ of $\lambda$. Using the Poincar\'e
duality isomorphism between homology and cohomology, $X_{\la}(F_\bull)$
defines a {\em Schubert class} $\s_{\la}=[X_{\la}(F_\bull)]$ in
$H^{2|\la|}(X,\Z)$. The algebraic group $GL_N(\C)$ acts transitively on $X$
and on the flags in $E$. The action of an element $g\in GL_N(\C)$ on the
variety $X_{\la}(F_\bull)$ is given by $g\cdot X_{\la}(F_\bull)=
X_{\la}(g\cdot F_\bull)$.
It follows that $\s_{\la}$ does not depend on the choice of flag
$F_\bull$ used to define $X_{\la}$. As all cohomology classes in these
lectures will occur in even degrees, we will adopt the convention that the
degree of a class $\alpha\in H^{2k}(X,\Z)$ is equal to $k$.

\medskip

We next review the classical facts about the cohomology of $X=G(m,N)$.

\medskip
\noin
{\bf 1)} The additive structure of $H^*(X,\Z)$ is given by 
\[
H^*(X,\Z)=\bigoplus_{\la\subset (n^m)} \Z\cdot \s_{\la},
\]
that is, $H^*(X,\Z)$ is a free abelian group with basis given by the
Schubert classes.

\begin{figure}[htpb]
\epsfxsize 35mm
\center{\mbox{\epsfbox{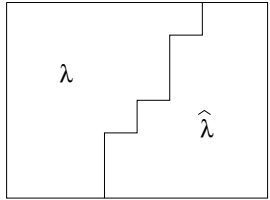}}}
\caption{Dual Young diagrams}
\label{duals}
\end{figure}

\medskip
\noin {\bf 2)} To describe the cup product in $H^*(X,\Z)$, we will use
Schubert's {Duality Theorem}. This states that for any $\la$ and $\mu$ with
$|\la|+|\mu|=mn$, we have $\s_{\la}\s_{\mu}= \delta_{\wh{\la}\mu}\cdot
[\mathrm{pt}]$, where $[\mathrm{pt}]=\s_{(n^m)}$ is the class of a point, and
$\wh{\la}$ is the dual partition to $\la$. The diagram of $\wh{\la}$ is the
complement of $\la$ in the rectangle $(n^m)$, rotated by $180^{\circ}$.  This
is illustrated in Figure \ref{duals}.

For general products, we have an equation
\[
\s_{\la}\, \s_{\mu}=\sum_{|\nu|=|\la|+|\mu|}c_{\la\mu}^{\nu} \,\s_{\nu}
\]
in $H^*(X,\Z)$, and the {\em structure constants} $c_{\la\mu}^{\nu}$ are
given by
\[
c_{\la\mu}^{\nu}=\int_X \s_{\la}\s_{\mu}\s_{\wh{\nu}} = \left<\s_{\la},
\s_{\mu},\s_{\wh{\nu}}\right>_0 = 
\# X_{\la}(F_\bull)\cap X_{\mu}(G_\bull) \cap X_{\wh{\nu}}(H_\bull),
\]
for general full flags $F_\bull$, $G_\bull$ and $H_\bull$ in $E$. Later, 
we will discuss a combinatorial formula for these structure constants.

\medskip
\noin
{\bf 3)} The classes $\s_1,\ldots,\s_n$ are called {\em special Schubert
classes}. Observe that there is a unique Schubert class in codimension
one: $H^2(X,\Z)=\Z\,\s_1$. If
\begin{equation}
\label{whitneysum}
0\ra S\ra E_X \ra Q \ra 0
\end{equation}
is the tautological short exact sequence of vector bundles over $X$, with
$E_X=X\times E$, then one can show that $\s_i$ is equal to the $i$th Chern
class $c_i(Q)$ of the quotient bundle $Q$, for $0\lequ i \lequ n$.

\begin{thm}[Pieri rule, \cite{Pi}] 
For $1\lequ p \lequ n$ we have $\s_{\la} \,\s_p=
\sum \s_{\mu}$, where the sum is over all $\mu\subset (n^m)$ obtained from
$\la$ by adding $p$ boxes, with no two in the same column.
\end{thm}

\begin{exa} Suppose $m=n=2$ and we consider the Grassmannian $X=G(2,4)$
of $2$-planes through the origin in $E=\C^4$. Note that $X$ may be identified
with the Grassmannian of all lines in projective $3$-space $P(E)\cong \bP^3$.
The list of Schubert classes for $X$ is
\[
\s_0=1, \ \  \s_1, \ \ \s_2, \ \ \s_{1,1}, \ \  \s_{2,1}, 
\ \ \s_{2,2}=[\mathrm{pt}].
\]
Observe that the indices of these classes are exactly the six partitions whose
diagrams fit inside a $2\times 2$ rectangle. Using the Pieri rule, we compute
that
\[
\s_1^2=\s_2+\s_{1,1}, \ 
\s_1^3 = 2\,\s_{2,1}, \ \s_1^4 = 2\,\s_{2,2}= 2\,[\mathrm{pt}].
\]
The last relation means that there are exactly $2$ points in the intersection
\[
X_1(F_\bull)\cap X_1(G_{\bull})\cap X_1(H_{\bull})\cap
X_1(I_\bull),
\]
for general flags $F_\bull$, $G_\bull$, $H_{\bull}$, and $I_{\bull}$. Since
e.g.\ $X_1(F_\bull)$ may be identified with the locus of lines in $P(E)$
meeting the fixed line $P(F_2)$, this proves the enumerative fact that there
are two lines in $\bP^3$ which meet four given lines in general position.
\end{exa}

\medskip
\noin {\bf 4)} Any Schubert class $\s_{\la}$ may be expressed as a polynomial
in the special classes in the following way. Let us agree here and in the
sequel that $\s_p=0$ if $p<0$ or $p>n$.

\begin{thm}[Giambelli formula, \cite{G}]
We have $\s_\la=\det(\s_{\la_i+j-i})_{1\lequ i,j\lequ m}$, that is, $\s_\la$
is equal to a Schur determinant in the special classes.
\end{thm}

\medskip
\noin
{\bf 5)} The ring $H^*(X,\Z)$ is presented as a quotient of the polynomial
ring $\Z[\s_1,\ldots,\s_n]$ by the relations 
\[
D_{m+1}=\cdots=D_N=0
\]
where $D_k=\det(\s_{1+j-i})_{1\lequ i,j\lequ k}$. To understand where 
these relations come from, note that the Whitney sum formula applied
to (\ref{whitneysum}) says that $c_t(S)c_t(Q)=1$, which implies, since
$\s_i=c_i(Q)$, that $D_k=(-1)^kc_k(S)=c_k(S^*)$. In particular, we 
see that $D_k$ vanishes for $k>m$, since $S^*$ is a vector bundle of 
rank $m$.

\subsection{Gromov-Witten invariants}
Our starting point is the aforementioned fact that the classical structure
constant $c_{\la\mu}^{\nu}$ in the cohomology of $X=G(m,N)$ can be realized as
a triple intersection number $\# X_{\la}(F_\bull)\cap X_{\mu}(G_\bull) \cap
X_{\wh{\nu}}(H_\bull)$ on $X$. The three-point, genus zero Gromov-Witten
invariants on $X$ extend these numbers to more general enumerative constants,
which are furthermore used to define the `small quantum cohomology ring' of
$X$.

A {\em rational map of degree} $d$ to $X$ is a morphism
$f\colon \bP^1\to X$ such that
$$\int_X f_*[\bP^1]\cdot \sigma_1=d,$$ 
i.e.\ $d$ is the number of points in
$f^{-1}(X_1(F_\bull))$ when $F_\bull$ is in general position. 

\begin{defn} 
Given a degree $d\gequ 0$ and partitions $\lambda$, $\mu$, and $\nu$ such that
$|\lambda| + |\mu| + |\nu| = mn + dN$, we define the Gromov-Witten invariant
$\langle \s_\lambda, \s_\mu, \s_\nu \rangle_d$ to be the
number of rational maps $f\colon \bP^1\to X$ of degree $d$ such that
$f(0)\in X_\lambda(F_\bull)$, $f(1)\in X_\mu(G_\bull)$,
and $f(\infty)\in X_\nu(H_\bull)$, for given flags $F_\bull$,
$G_\bull$, and $H_\bull$ in general position. 
\end{defn}
\noin
We shall show later that $\langle \s_\lambda, \s_\mu, \s_\nu \rangle_d$
is a well-defined, finite integer. Notice that for 
the degree zero invariants, we have
\[
\langle \s_\lambda, \s_\mu, \s_\nu \rangle_0= 
\int_X \s_{\lambda} \s_\mu \s_\nu=
\# X_{\la}(F_\bull)\cap X_{\mu}(G_\bull) \cap
X_{\nu}(H_\bull),
\]
as a morphism of degree zero is just a constant map to $X$.

\medskip
\noin {\bf Key example}.  Consider the Grassmannian
$G(d,2d)$ for any $d\gequ 0$. We say that two points $U$, $V$ of $G(d,2d)$ are
in general position if the intersection $U\cap V$ of the corresponding
subspaces is the zero subspace.

\begin{prop}[\cite{BKT1}]
\label{3pts:typea}
Let $U$, $V$, and $W$ be three points of $\zZ=G(d,2d)$ which are pairwise in
general position.  Then there is a unique morphism $f\colon \bP^1\ra \zZ$ of
degree $d$ such that $f(0)=U$, $f(1)=V$, and $f(\infty)=W$. In particular, the
Gromov-Witten invariant which counts degree $d$ maps to $\zZ$ through three
general points is equal to $1$.
\end{prop}
\begin{proof}
Let $U$, $V$, and $W$ be given, in pairwise general position.  Choose a basis
$(v_1, \dots, v_d)$ of $V$. Then we can construct a morphism $f \colon \bP^1
\to \zZ$ of degree $d$ such that $f(0) = U$, $f(1) = V$, and $f(\infty) = W$
as follows.  For each $i$ with $1 \lequ i \lequ d$, we let $u_i$ and $w_i$ be
the projections of $v_i$ onto $U$ and $W$, respectively.  If $(s \cc t)$ are
the homogeneous coordinates on $\bP^1$, then the morphism 
\[
f(s \cc t) = \Span\{ s u_1 + t w_1, \dots, s u_d + t w_d \}
\] 
satisfies the required conditions.  Observe that $f$ does not depend on the
chosen basis for $V$. Indeed, if $v_i'=\sum a_{ij}v_j$, then $u_i'=\sum
a_{ij}u_j$, $w_i'=\sum a_{ij} w_j$ and one checks easily that
\[
 \Span \{ s u_1 + t w_1, \dots, s u_d + t w_d \} = 
 \Span \{ s u'_1 + t w'_1, \dots, s u'_d + t w'_d \}.
\]

\medskip
\noin
{\bf Exercise.} Show that the map $f$ is an embedding of 
$\bP^1$ into $\zZ$ such that $f(p_1)$ and $f(p_2)$ are in general position,
for all points $p_1$, $p_2$ in $\bP^1$ with $p_1\neq p_2$. Show also that
$f$ has degree $d$.

\medskip
Next, suppose that $f \colon \bP^1 \to \zZ$ is any morphism of degree $d$
which sends $0$, $1$, $\infty$ to $U$, $V$, $W$, respectively.  Let $S \subset
\C^{2d} \otimes \cO_{\zZ}$ be the tautological rank $d$ vector bundle over
$\zZ$, and consider the pullback $f^*S\ra \bP^1$.  The morphism $f \colon
\bP^1 \to G(d,2d)$ is determined by the inclusion of $f^*S$ in
$\C^{2d}\otimes\cO_{\bP^1}$, i.e., a point $p\in \bP^1$ is mapped by $f$ to
the fiber over $p$ of the image of this inclusion.

Every vector bundle over $\bP^1$ splits as a direct sum of line bundles, so
$f^*S\cong \oplus_{i=1}^d\cO(a_i)$. Each $\cO(a_i)$ is a subbundle of a
trivial bundle, hence $a_i\lequ 0$, and $\sum a_i=-d$ as $f$ has degree
$d$. We deduce that $a_i=-1$ for each $i$, since otherwise $f^*S$ would have a
trivial summand, and this contradicts the general position hypothesis. It
follows that we can write $f(s \cc t) = \Span \{ s u_1 + t w_1, \dots, s u_d +
t w_d \}$ for suitable vectors $u_i, w_i \in \C^{2d}$, which depend on the
chosen identification of $f^*S$ with $\oplus_{i=1}^d \cO(-1)$. We conclude
that $f$ is the map constructed as above from the basis $(v_1,\dots,v_d)$,
where $v_i=u_i+w_i$.
\end{proof}

\medskip

We now introduce the key definition upon which the subsequent analysis depends.

\begin{defn}[\cite{Buch}] For any morphism $f\colon \bP^1\to G(m,N)$, 
define the {\em kernel\/} of $f$
to be the intersection of all the subspaces $V \subset E$ corresponding to
image points of $f$.  Similarly, the {\em span\/} of $f$ is the linear span of
these subspaces.
\[ \Ker(f) = \bigcap_{p \in \bP^1} f(p); \hspace{.5cm}
   \Span(f) = \sum_{p \in \bP^1} f(p).
\]
\end{defn}
\noin
Note that for each $f\colon \bP^1\to X$, we have 
$\Ker(f)\subset \Span(f)\subset E$.

\begin{lemma}[\cite{Buch}] 
\label{buchlemma}
If $f\colon \bP^1\to G(m,N)$ is a morphism of 
degree $d$, then $\dim \Ker(f)\gequ m-d$ and $\dim\Span(f)\lequ m+d$.
\end{lemma}
\begin{proof}
Let $S\ra X$ be the rank $m$ tautological bundle over $X=G(m,N)$.  Given any
morphism $f\colon \bP^1\to X$ of degree $d$, we have that $f^*S\cong
\oplus_{i=1}^m\cO(a_i)$, where $a_i\lequ 0$ and $\sum a_i=-d$. Moreover, the
map $f$ is induced by the inclusion $f^*S\subset
E\otimes\cO_{\bP^1}$. There are at least $m-d$ zeroes among the integers
$a_i$, hence $f^*S$ contains a trivial summand of rank at least $m-d$. But
this corresponds to a fixed subspace of $E$ of the same dimension which is
contained in $\Ker(f)$, hence $\dim \Ker(f)\gequ m-d$.

Similar reasoning shows that if $Q \ra X$ is the rank
$n$ universal 
quotient bundle over $X$, then $f^*Q$ has a trivial summand of rank at
least $n-d$. It follows that the image of the map 
$f^*S\ra E\otimes\cO_{\bP^1}$ factors through a subspace of $E$
of codimension at least $n-d$, and hence of dimension at most $m+d$.
\end{proof}

In the next lecture, we will see that for those maps $f$ which are 
counted by a degree $d$ Gromov-Witten invariant for $X$, we have
$\dim\Ker(f)=m-d$ and $\dim\Span(f)=m+d$. In fact, it will turn out
that the pair $(\Ker(f), \Span(f))$ determines $f$ completely!

\newpage

\section{Lecture Two}

\subsection{The main theorem}
Given integers $a$ and $b$, we let $F(a,b;E)=F(a,b; N)$ denote the two-step
flag variety parametrizing pairs of subspaces $(A,B)$ with $A \subset B
\subset E$, $\dim A = a$ and $\dim B = b$. We agree that $F(a,b;N)$ is empty
unless $0\lequ a\lequ b\lequ N$; when the latter condition holds then
$F(a,b;N)$ is a projective complex manifold of dimension $(N-b)b+(b-a)a$.  For
any non-negative integer $d$ we set $Y_d = F(m-d, m+d; E)$; this will be the
parameter space of the pairs $(\Ker(f),\Span(f))$ for the relevant morphisms
$f\colon \bP^1\ra X$.  Our main theorem will be used to identify Gromov-Witten
invariants on $X=G(m,E)$ with classical triple intersection numbers on the
flag varieties $Y_d$.

To any subvariety $\W \subset X$ we associate the subvariety
$W^{(d)}$ in $Y_d$ defined by
\begin{equation}
\label{assocsubvar}
\W^{(d)} = \{\,(A,B) \in Y_d \, \mid\, \exists~V \in \W : A \subset V
   \subset B \,\} \,.
\end{equation}
Let $F(m-d,m,m+d;E)$ denote the variety of three-step flags in $E$ of 
dimensions $m-d$, $m$, and $m+d$. There are natural projection maps
\[
\pi_1 \colon F(m-d,m,m+d;E) \to X  \ \ \mathrm{ and } \ \ \pi_2\colon
F(m-d,m,m+d;E) \to Y_d.
\]
We then have $\W^{(d)} =
\pi_2(\pi_1^{-1}(\W))$. Moreover, as the maps $\pi_i$ are $GL_N$-equivariant,
if $\W = X_\lambda(F_\bull)$ is a Schubert variety in $X$, then $\W^{(d)} =
X_\lambda^{(d)}(F_\bull)$ is a Schubert variety in $Y_d$. We will describe this
Schubert variety in more detail after we prove the main theorem.  

\medskip
\noin
{\bf Remarks.} 1) One computes that $\dim Y_d=mn+dN-3d^2$.

\medskip
\noin
2) Since the fibers of $\pi_2$ are isomorphic to $G(d,2d)$, the codimension of
$X_\lambda^{(d)}(F_\bull)$ in $Y_d$ is at least $|\lambda|-d^2$. 
 
\begin{thm}[\cite{BKT1}] \label{thm:typea}
Let $\lambda$, $\mu$, and $\nu$ be partitions and $d$ be an integer
such that $|\lambda| + |\mu| + |\nu| = mn + dN$,
and let $F_\bull$, $G_\bull$, and $H_\bull$ be complete flags of
$E = \C^N$ in general position.  Then the
map $f \mapsto (\Ker(f),\Span(f))$ gives a bijection of the set of
rational maps $f\colon \bP^1 \to G(m,N)$ of degree $d$ satisfying
$f(0)\in X_\lambda(F_\bull)$, $f(1)\in X_\mu(G_\bull)$, and
$f(\infty)\in X_\nu(H_\bull)$, with the set of points 
in the intersection
$X^{(d)}_\lambda(F_\bull) \cap X^{(d)}_\mu(G_\bull) \cap
X^{(d)}_\nu(H_\bull)$ in $Y_d=F(m-d,m+d;N)$.
\end{thm}

It follows from Theorem \ref{thm:typea} that 
we can express any 
Gromov-Witten invariant of degree $d$ on $G(m,N)$ as
a classical intersection number on $Y_d$.  Let
$[X_\lambda^{(d)}]$ denote the cohomology class of
$X^{(d)}_\lambda(F_\bull)$ in $H^*(Y_d,\Z)$.

\begin{cor}
\label{cor:typea}
Let $\lambda$, $\mu$, and $\nu$ be partitions and $d\gequ 0$  an integer
such that $|\lambda| + |\mu| + |\nu| = mn + dN$. We then have
\[
\langle \s_\lambda, \s_\mu, \s_\nu \rangle_d=
\int_{F(m-d,m+d;N)} 
[X^{(d)}_{\lambda}]\cdot [X^{(d)}_{\m}]\cdot [X^{(d)}_{\nu}].
\]
\end{cor}

\begin{proof}[Proof of Theorem \ref{thm:typea}] 
Let $f\colon \bP^1\to X$ be a rational map as in
the statement of the theorem.

\medskip
\noin
{\bf Claim 1.} We have $d\lequ \min(m,n)$,
$\dim\Ker(f)=m-d$ and $\dim\Span(f)=m+d$.

\medskip
Indeed, let $a=\dim\Ker(f)$ and $b=\dim\Span(f)$.
In the two-step flag variety $Y'=F(a,b;E)$ there are 
associated Schubert varieties $X'_\lambda(F_{\bull})$, $X'_\mu(G_{\bull})$, 
and $X'_\nu(H_{\bull})$, defined as in (\ref{assocsubvar}).
Writing $e_1=m-a$ and $e_2=b-m$, we see that the codimension of
$X'_\lambda(F_{\bull})$ in $Y'$ is at least  $|\lambda|-e_1e_2$, 
and similar inequalities hold with $\mu$ and $\nu$ in place of $\lambda$.
Since $(\Ker(f),\Span(f))\in X'_\lambda(F_{\bull})\cap X'_\mu(G_{\bull})\cap
X'_\nu(H_{\bull})$ and the three flags $F_\bull$, $G_\bull$ and $H_\bull$
are in general position, we obtain
\[
mn+dN-3e_1e_2\lequ \dim F(a,b;E)=(N-b)(m+e_2)+(e_1+e_2)a,
\]
and hence, by a short computation,
\begin{equation}
\label{ineq1}
dN\lequ 2e_1e_2+e_2(N-b)+ae_1.
\end{equation}
Lemma \ref{buchlemma} says that $e_1\lequ d$ and $e_2\lequ d$, and therefore
that the right-hand side of (\ref{ineq1}) is at most $2e_1e_2+d(N-b+a)$. Since
$b-a=e_1+e_2$, it follows that
\[
(e_1+e_2)^2\lequ 2d(e_1+e_2)\lequ 4e_1e_2,
\] 
and hence $e_1=e_2=d$. This proves Claim 1.

\medskip

Let $\CC$ denote the 
set of rational maps in the statement of the theorem,
and set $\II = X^{(d)}_\lambda(F_\bull) \cap
X^{(d)}_\mu(G_\bull) \cap X^{(d)}_\nu(H_\bull)$.
If $f \in \CC$ then Claim 1 shows that 
$(\Ker(f),\Span(f)) \in \II$. We next describe the
inverse of the resulting map $\CC\ra\II$.

Given $(A,B) \in \II$, we let $\zZ = G(d,B/A) \subset X$ be the set of
$m$-dimensional subspaces of $E$ between $A$ and $B$.  Observe that $\zZ\cong
G(d,2d)$, and that $X_\lambda(F_\bull) \cap \zZ$, $X_\mu(G_\bull) \cap \zZ$,
and $X_\nu(H_\bull) \cap \zZ$ are non-empty Schubert varieties in
$\zZ$. (Indeed, e.g.\ $X_\lambda(F_\bull) \cap \zZ$ is defined by the attitude
of $V/A$ with respect to the flag $\wt{F}_\bull$ in $B/A$ with
$\wt{F}_i=((F_i+A)\cap B)/A$ for each $i$.) We assert that each of
$X_\lambda(F_\bull) \cap \zZ$, $X_\mu(G_\bull) \cap \zZ$, and $X_\nu(H_\bull)
\cap \zZ$ must be a single point, and that these three points are subspaces 
of $B/A$ in pairwise
general position.  Proposition \ref{3pts:typea} then provides the unique
$f\colon \bP^1\to X$ in $\CC$ with $\Ker(f)=A$ and $\Span(f)=B$.

\medskip
\noin
{\bf Claim 2.} Let $U$, $V$, and $W$ be three points in $\zZ$, one from each 
intersection. Then the subspaces $U$, $V$ and $W$ are in pairwise general
position.

\medskip
Assuming this claim, we can finish the proof as follows.  Observe that any
positive dimensional Schubert variety in $\zZ$ must contain a point $U'$ which
meets $U$ non-trivially, and similarly for $V$ and $W$. Indeed, on $G(d,2d)$,
the locus of $d$-dimensional subspaces $\Sigma$ with $\Sigma\cap U\neq\{0\}$
is, up to a general translate, the unique Schubert variety in codimension
$1$. It follows that this locus must meet any other Schubert variety
non-trivially, unless the latter is zero dimensional, in other words, a
point. Therefore Claim 2 implies that $X_\lambda(F_\bull) \cap \zZ$,
$X_\mu(G_\bull) \cap \zZ$, and $X_\nu(H_\bull) \cap \zZ$ are three points on
$\zZ$ in pairwise general position.

To prove Claim 2, we again use a dimension counting argument to show that if
the three reference flags are chosen generically, no two subspaces among $U$,
$V$, $W$ can have non-trivial intersection.  Consider the three-step flag
variety $Y''=F(m-d,m-d+1,m+d; E)$ and the projection $\pi\colon Y''\ra Y_d$.
Note that $\dim Y''=\dim Y_d+2d-1$, as $Y''$ is a $\bP^{2d-1}$-bundle over
$Y_d$.  

To each subvariety $\W\subset G(m,E)$ we associate $\W''\subset Y''$
defined by
\[
\W'' = \{\,(A,A',B) \in Y'' \,\mid\, \exists~V \in \W : A' \subset V
   \subset B \,\} \,.
\]
We find that the codimension of $X''_\mu(G_{\bull})$ in $Y''$ is at least 
$|\mu|-d^2+d$, and similarly for $X''_\nu(H_{\bull})$. 
Since the three flags are in general position, and
$\pi^{-1}(X^{(d)}_{\la}(F_{\bull}))$ has codimension 
at least $|\la|-d^2$ in $Y''$,
we must have
\[
\pi^{-1}(X^{(d)}_{\la}(F_{\bull})) 
\cap X''_\mu(G_{\bull}) \cap X''_\nu(H_{\bull})=\emptyset,
\]
and the same is true for the other two analogous triple intersections.
This completes the proof of Claim 2, and of the theorem.
\end{proof}

It is worth pointing out that we may rephrase Theorem 
\ref{thm:typea}
using rational curves in $X$, instead of rational maps to $X$.
For this, recall from the Exercise given in the first lecture 
that every rational map $f$ that is counted
in Theorem \ref{thm:typea} is an
embedding of $\bP^1$ into $X$ of degree equal to the degree
of the curve $\mathrm{Im}(f)$. Moreover, the bijection of the theorem
shows that all of these maps have different images.

\subsection{Parametrizations of Schubert varieties}
We now describe an alternative way to parametrize the Schubert 
varieties on $G(m,N)$, by replacing each partition $\lambda$
by a $01$-string $I(\lambda)$ of length $N$, with $m$ zeroes.
Begin
by drawing the Young diagram of the partition $\lambda$ in the
upper-left corner of an $m \times n$ rectangle. We then put a 
label on each step of the path from the
lower-left to the upper-right corner of this rectangle which follows
the border of $\lambda$. Each vertical step is labeled ``0'', while
the remaining $n$ horizontal steps are labeled ``1''. 
The string $I(\lambda)$ then consists of these labels in
lower-left to upper-right order.  

\begin{exa} 
\label{grex}
On the Grassmannian $G(4,9)$, the $01$-string of the 
partition $\la=(4,4,3,1)$ is $I(\la)=101101001$. This is
illustrated below.
$$\includegraphics[scale=.9]{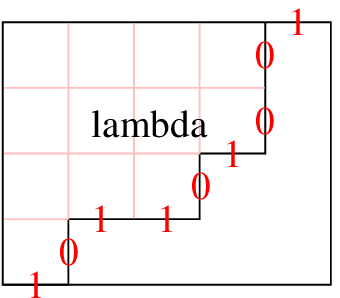}$$
\end{exa}

Alternatively, each partition $\la\subset (n^m)$ corresponds to a 
Grassmannian permutation $w_{\la}$ in the symmetric group $S_N$, which
is a minimal length representative in the coset space $S_N/(S_m\times
S_n)$. The element $w=w_{\la}$ is such that the positions of the ``0''s
(respectively, the ``1''s) in the $01$-string $I(\la)$ are given 
by $w(1),\ldots,w(m)$ (respectively, by $w(m+1),\ldots,w(N)$). In Example
\ref{grex}, we have $w_{\la}=257813469\in S_9$. 

\medskip

In a similar fashion, the Schubert varieties on the two-step flag variety
$F(a,b;N)$ are parametrized by permutations $w\in S_N$ with $w(i)<w(i+1)$ for
$i\notin \{a,b\}$.  For each such permutation $w$ and fixed full flag
$F_{\bull}$ in $E$, the Schubert variety $X_w(F_{\bull})\subset F(a,b;N)$ is
defined as the locus of flags $A\subset B \subset E$ such that
\[
\dim(A\cap F_i) \gequ \#\{p\lequ a\ |\ w(p) > N-i\} \ \ \mathrm{and} \ \ 
\dim(B\cap F_i) \gequ \#\{p\lequ b\ |\ w(p)> N-i\}
\]
for each $i$. The codimension of $X_w(F_{\bull})$ in $F(a,b;N)$ is equal to
the {\em length} $\ell(w)$ of the permutation $w$. Furthermore, these indexing
permutations $w$ correspond to $012$-strings $J(w)$ of length $N$ with $a$
``0''s and $b-a$ ``1''s. The positions of the ``0''s (respectively, the
``1''s) in $J(w)$ are recorded by $w(1),\ldots,w(a)$ (respectively, by
$w(a+1),\ldots,w(b)$).

\medskip
Finally, we describe the $012$-string $J^d(\lambda)$ 
associated to the modified 
Schubert variety $X^{(d)}_\la(F_\bull)$ in $Y_d = F(m-d,m+d;N)$.
This string is obtained by first multiplying each number in 
the $01$-string $I(\la)$ by $2$, to get a $02$-string $2 I(\la)$. 
We then get the $012$-string $J^d(\lambda)$ by
changing the first $d$ ``2''s and the last $d$ ``0''s of 
$2 I(\la)$ to ``1''s. Taking $d=2$ in Example \ref{grex}, 
we get $J^2(4,4,3,1)=101202112$, which corresponds to a Schubert
variety in $F(2,6;9)$. 

\begin{cor}[\cite{Y}]
\label{cor2:typea}
Let $\lambda$, $\mu$, and $\nu$ be partitions and $d\gequ 0$ be 
such that $|\lambda| + |\mu| + |\nu| = mn + dN$. If any of 
$\lambda_d$, $\mu_d$, and $\nu_d$ is less than $d$, then
$\langle \s_\lambda, \s_\mu, \s_\nu \rangle_d=0$.
\end{cor}
\begin{proof}
By computing the length of the permutation corresponding to $J^d(\la)$, one
checks easily that when $\lambda_d<d$, the codimension of
$X^{(d)}_\lambda(F_\bull)$ in $Y_d$ is strictly greater than $|\lambda|-d^2$.
Therefore, when any of $\lambda_d$, $\mu_d$, or $\nu_d$ is less than $d$, the
sum of the codimensions of the three Schubert varieties
$X^{(d)}_{\la}(F_\bull)$, $X^{(d)}_{\mu}(G_\bull)$, and
$X^{(d)}_{\nu}(H_\bull)$ which appear in the statement of Theorem
\ref{thm:typea} is strictly greater than the dimension of $Y_d=F(m-d,m+d;N)$.
We deduce that
\[
\langle \s_\lambda, \s_\mu, \s_\nu \rangle_d=
\int_{Y_d}[X_{\la}^{(d)}]\cdot[X_{\mu}^{(d)}]\cdot [X_{\nu}^{(d)}]=0.
\]
\end{proof}

\newpage

\section{Lecture Three} 

\subsection{Classical and quantum Littlewood-Richardson rules}
The problem we turn to now is that of finding a positive combinatorial formula
for the Gromov-Witten invariants $\langle \s_\lambda, \s_\mu, \s_\nu
\rangle_d$. For the classical structure constants (the case $d=0$), this
problem was solved in the 1930's by Littlewood and Richardson, although
complete proofs only appeared in the early 1970's. In the past few years,
there has been a resurgence of interest in this question (see, for example,
\cite{Fu}), which has led to a new formulation of the rule in terms of
`puzzles' (due to Knutson, Tao, and Woodward).

Define a {\em puzzle\/} to be a triangle decomposed into {\em puzzle
pieces\/} of the three types displayed below. 
$$\pzp{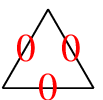} \pzs \pzs \pzp{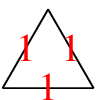}
\pzs \pzs \raisebox{-2.5mm}{\pzp{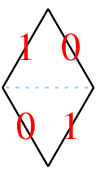}} $$
A puzzle piece may be rotated but not reflected when used in a puzzle.
Furthermore, the common edges of two puzzle pieces next
to each other must have the same labels. Recall from the last lecture
that a Schubert class $\s_{\la}$ in $H^*(G(m,N),\Z)$ may also 
be indexed by a $01$-string $I(\la)$ with $m$ ``0''s and $n$ ``1''s. 

\begin{thm}[\cite{KTW}] \label{thm:ktw}
For any three Schubert classes $\s_{\la}$, $\s_{\mu}$, and $\s_{\nu}$ 
in the cohomology of $X=G(m,N)$, the integral $\int_X \s_{\la}\s_{\mu}
\s_{\nu}$ is equal to the number of puzzles such that  
$I(\la)$, $I(\mu)$, and
$I(\nu)$ are the labels on the north-west, north-east, and
south sides when read in clockwise order.
\end{thm}
\noin The formula in Theorem \ref{thm:ktw} is bijectively equivalent to the
classical Littlewood-Richardson rule, which describes the same numbers as the
cardinality of a certain set of Young tableaux (see \cite[\S 4.1]{V}).

\begin{exa} In the projective plane $\bP^2=G(1,3)$, two general lines intersect
in a single point. This corresponds to the structure constant
\[
\langle \s_1, \s_1, \s_0 \rangle_0=
\int_{\bP^2} \s_1^2= 1.
\]
The figure below displays the unique puzzle with the corresponding three
strings $101$, $101$, and $011$ on its north-west, north-east, and south sides.
$$\includegraphics[scale=.9]{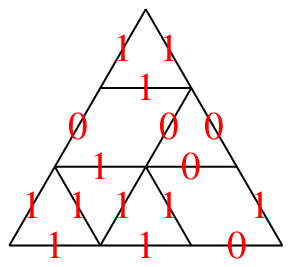}$$ We suggest that the reader works
out the puzzle which corresponds to the intersection
$\s_k\s_{\ell}=\s_{k+\ell}$ on $\bP^n$, for $k+\ell\lequ n$.
\end{exa}


It is certainly tempting to try to generalize Theorem \ref{thm:ktw} to
a result that would hold for the flag variety $SL_N/B$. This time the
three sides of the puzzle would be labeled by permutations, and one has
to specify the correct set of puzzle pieces to make the rule work.
In the fall of 1999, Knutson proposed such a general conjecture for
the Schubert structure constants on all partial flag varieties, which 
specialized to Theorem \ref{thm:ktw} in the Grassmannian case. However,
he soon discovered counterexamples to this conjecture (in 
fact, it fails for the three-step flag variety $F(1,3,4;5)$). 

Motivated by Theorem \ref{thm:typea}, Buch, Kresch and the author were
especially interested in a combinatorial rule for the structure constants on
two-step flag varieties. Surprisingly, there is extensive computer evidence
which suggests that Knutson's conjecture is true in this special case. Recall
from the last lecture that the Schubert classes on two-step flag varieties are
indexed by the $012$-strings $J(w)$, for permutations $w\in S_N$. 
In this setting we have the following six different
types of puzzle pieces.

$$\pzp{gwistri0} \pzs \pzp{gwistri1} \pzs \pzp{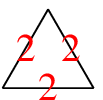} \pzs \pzp{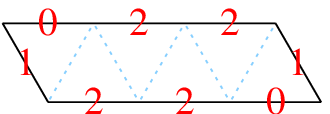}
\pzs \raisebox{-2.5mm}{\pzp{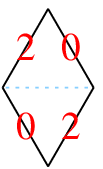}} \pzs \pzp{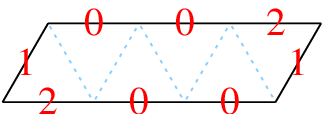}$$

\begin{figure}
$$\includegraphics[scale=.9]{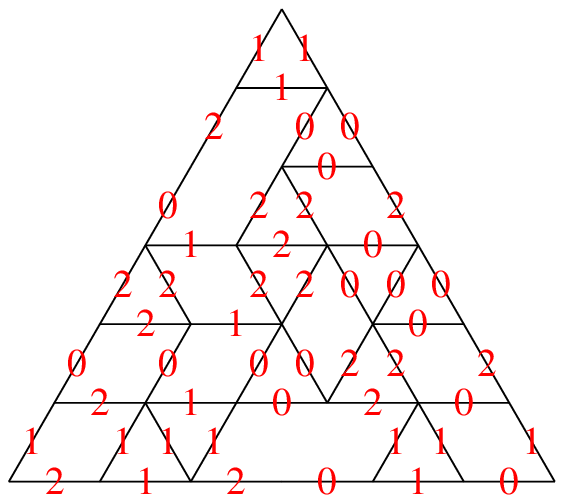}
\includegraphics[scale=.9]{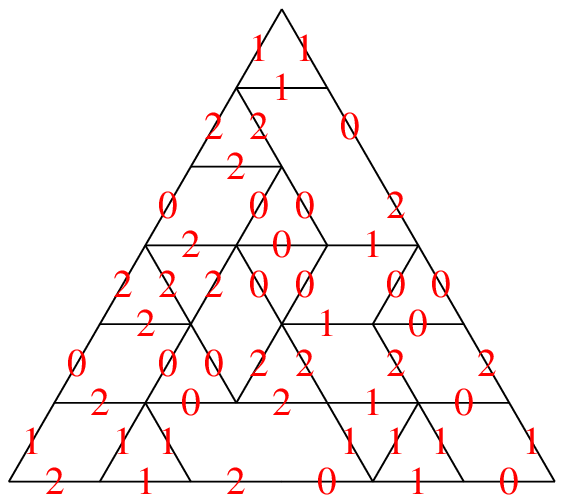}$$
\caption{Two puzzles with the same boundary labels}
\label{puzzleexa}
\end{figure}

The length of the fourth
and of the sixth piece above may vary.  The fourth piece can
have any number of ``2''s (including none) to the right of the ``0''
on the top edge and equally many to the left of the ``0'' on the
bottom edge.  Similarly the sixth piece can have an arbitrary number
of ``0''s on the top and bottom edges.
Again each puzzle piece may be rotated but not reflected.
Figure \ref{puzzleexa} shows two examples of such puzzles.

We can now state Knutson's conjecture in the case of two-step flag
varieties. This conjecture has been verified by computer for all two-step flag
varieties $F(a,b;N)$ for which $N \lequ 16$.

\begin{conj}[Knutson] \label{conj:knutson}
For any three Schubert varieties $X_u$, $X_v$, and $X_w$ 
in the flag variety $F(a,b;N)$, the integral $\int_{F(a,b;N)}
[X_u]\cdot [X_v]\cdot [X_w]$ is equal to the number of puzzles such that  
$J(u)$, $J(v)$, and
$J(w)$ are the labels on the north-west, north-east, and
south sides when read in clockwise order.
\end{conj}

By combining Theorem \ref{thm:typea} with Conjecture \ref{conj:knutson},
we arrive at a conjectural `quantum Littlewood-Richardson rule' for the
Gromov-Witten invariants $\langle\s_\lambda,
    \s_\mu, \s_\nu\rangle_d$. This time we use the 
$012$-string $J^d(\lambda)$ associated to the 
Schubert variety $X^{(d)}_\la(F_\bull)$ in $F(m-d,m+d;N)$.

\begin{conj}[\cite{BKT1}]
\label{qlrconj}
  For partitions $\lambda,\mu,\nu$ such that $|\lambda|+|\mu|+|\nu| =
  mn + dN$ the Gromov-Witten invariant $\langle\s_\lambda,
    \s_\mu, \s_\nu\rangle_d$ is equal to the number of puzzles such
  that $J^d(\lambda)$, $J^d(\mu)$, and $J^d(\nu)$ are the labels on the
  north-west, north-east, and south sides when read in clockwise
  order.
\end{conj}

The verified cases of Conjecture~\ref{conj:knutson} imply that 
Conjecture~\ref{qlrconj} holds for all Grassmannians
$G(m,N)$ for which $N \lequ 16$.  It has also been proved in some
special cases including when $\lambda$ has
length at most $2$ or when $m$ is at most $3$.  
 
\begin{exa}
On the Grassmannian $G(3,6)$, the Gromov-Witten invariant 
\[
\langle \s_{3,2,1},\, \s_{3,2,1}, \, \s_{2,1}\rangle_1
\] 
is equal to 
$2$. We have $J^1(3,2,1)= 102021$ and $J^1(2,1)=010212$.  
Figure \ref{puzzleexa} displays the 
two puzzles with the labels $J^1(3,2,1)$, $J^1(3,2,1)$, and 
$J^1(2,1)$ on their sides.
\end{exa}

\subsection{Quantum cohomology of $G(m,N)$} As was alluded to 
in the first lecture and also by the phrase `quantum Littlewood-Richardson
rule', the above Gromov-Witten invariants are the structure constants in a
deformation of the cohomology ring of $X=G(m,N)$. This (small) {\em quantum
cohomology ring} $QH^*(X)$ was introduced by string theorists, and is a
$\Z[q]$-algebra which is isomorphic to $H^*(X,\Z)\otimes_{\Z}\Z[q]$ as a
module over $\Z[q]$. Here $q$ is a formal variable of degree $N=m+n$. The ring
structure on $QH^*(X)$ is determined by the relation
\begin{equation}
\label{qprod}
\s_{\la}\cdot \s_{\m} = 
\sum \langle \s_{\la}, \s_{\m}, \s_{\wh{\n}} \rangle_d\,\s_{\n}\, q^d,
\end{equation}
the sum over $d\gequ 0$ and partitions $\n$ with
$|\n|=|\la|+|\m|-dN$. Note that the terms corresponding to $d=0$ just
give the classical cup product in $H^*(X,\Z)$. We will need to use
the hard fact that equation (\ref{qprod}) defines an {\em 
associative} product, which turns $QH^*(X)$ into a commutative ring with
unit. The reader can find a proof of 
this basic result in the expository paper \cite{FP}.

We now prove, following \cite{Buch}, analogues of the basic structure theorems
about $H^*(X,\Z)$ for the quantum cohomology ring $QH^*(X)$. For any Young
diagram $\lambda\subset (n^m)$, let $\ov{\la}$ denote the diagram obtained by
removing the leftmost $d$ columns of $\lambda$. In terms of partitions, we
have $\ov{\la}_i=\max\{\la_i-d,0\}$. For any Schubert variety
$X_{\la}(F_\bull)$ in $G(m,E)$, we consider an associated Schubert variety
$X_{\ov{\la}}(F_\bull)$ in $G(m+d,E)$. It is easy to see that if
$\pi:F(m-d,m+d;E)\ra G(m+d,E)$ is the projection map, then
$\pi(X_{\la}^{(d)}(F_\bull))=X_{\ov{\la}}(F_\bull)$.

\begin{cor} 
\label{cor1}
If $\langle \s_{\la}, \s_{\m}, \s_{\nu} \rangle_d\neq 0$, then 
$[X_{\ov{\la}}]\cdot [X_{\ov{\mu}}]\cdot [X_{\ov{\nu}}] \neq 0$
in $H^*(G(m+d,E),\Z)$. 
\end{cor}

\begin{cor} 
\label{cor2}
If $\langle \s_{\la}, \s_{\m}, \s_{\nu} \rangle_d\neq 0$ and 
$\ell(\la)+\ell(\mu)\lequ m$, then $d=0$.
\end{cor}
\begin{proof}
We know a priori that $|\la|+|\mu|+|\nu|=mn+dN$. The assumption on
the lengths of $\la$ and $\mu$ implies that
\[
|\ov{\la}|+|\ov{\mu}|+|\ov{\nu}|\gequ 
|\la|+|\mu|+|\nu|-2md = \dim G(m+d,E)+d^2.
\]
By Corollary \ref{cor1}, we must have $d=0$.
\end{proof}

Corollary \ref{cor2} implies that if $\ell(\la)+\ell(\mu)\lequ m$, then 
\[
\s_{\la}\cdot \s_{\mu}=
\sum_{d,\nu}\langle \s_{\la}, \s_{\m}, \s_{\wh{\n}} \rangle_d\,\s_{\n}\, q^d
=\sum_{|\nu|=|\la|+|\mu|}
\langle \s_{\la}, \s_{\m}, \s_{\wh{\n}} \rangle_0\,\s_{\n}, 
\]
that is, there are no quantum correction terms in the product 
$\s_{\la}\cdot \s_{\mu}$. 

\begin{thm}[Quantum Giambelli, \cite{Ber}] 
\label{qgiamba}
We have 
$\s_{\la}=\det(\s_{\la_i+j-i})_{1\lequ i,j\lequ m}$ in $QH^*(X)$. That
is, the classical Giambelli and quantum Giambelli formulas coincide for
$G(m,N)$.
\end{thm}
\begin{proof}
Define a linear map $\phi:H^*(X,\Z)\ra QH^*(X)$ by $\phi([X_{\la}])=
\s_{\la}$. It follows from Corollary \ref{cor2} that
$\s_p\cdot\s_{\mu}=\phi([X_p] [X_\mu])$ whenever $\ell(\mu)\lequ m-1$.
Using the classical Pieri rule and induction, we see that
\[
\s_{p_1}\cdots \s_{p_m}= \phi([X_{p_1}]\cdots [X_{p_m}]),
\]
for any $m$ special Schubert classes $\s_{p_1},\ldots ,\s_{p_m}$. This
implies that
\[
\det(\s_{\la_i+j-i})_{1\lequ i,j\lequ m}=
\phi(\det([X_{\la_i+j-i}])_{1\lequ i,j \lequ m})=
\phi([X_{\la}])= \s_{\la}.
\]
\end{proof}

\begin{thm}[Quantum Pieri, \cite{Ber}] 
\label{qpieria}
For $1\lequ p \lequ n$, we have
\begin{equation}
\label{pr}
\s_{\la}\cdot\s_p=\sum_{\mu}\s_{\mu}+q\sum_{\nu}\s_{\nu},
\end{equation}
where the first sum is over diagrams $\mu$ obtained from $\la$ by adding 
$p$ boxes, no two in the same column, and the second sum is over all $\nu$
obtained from $\la$ by removing $N-p$ boxes from the `rim' of $\lambda$, 
at least one from each row.
\end{thm}
\noin
Here the `rim' of a diagram $\lambda$ is the rim hook (or `border strip')
contained in $\la$ whose south-east border follows 
the path we used earlier to define the 
$01$-string corresponding to $\la$.

\begin{exa} For the Grassmannian $G(3,6)$, we have
\[
\s_{3,2,1} \cdot \s_2 = \s_{3,3,2}+q(\s_2+\s_{1,1})
\]
in $QH^*(G(3,6))$. The rule for obtaining
the two $q$-terms is illustrated below.
\skipline
\epsfxsize 65mm
\center{\mbox{\epsfbox{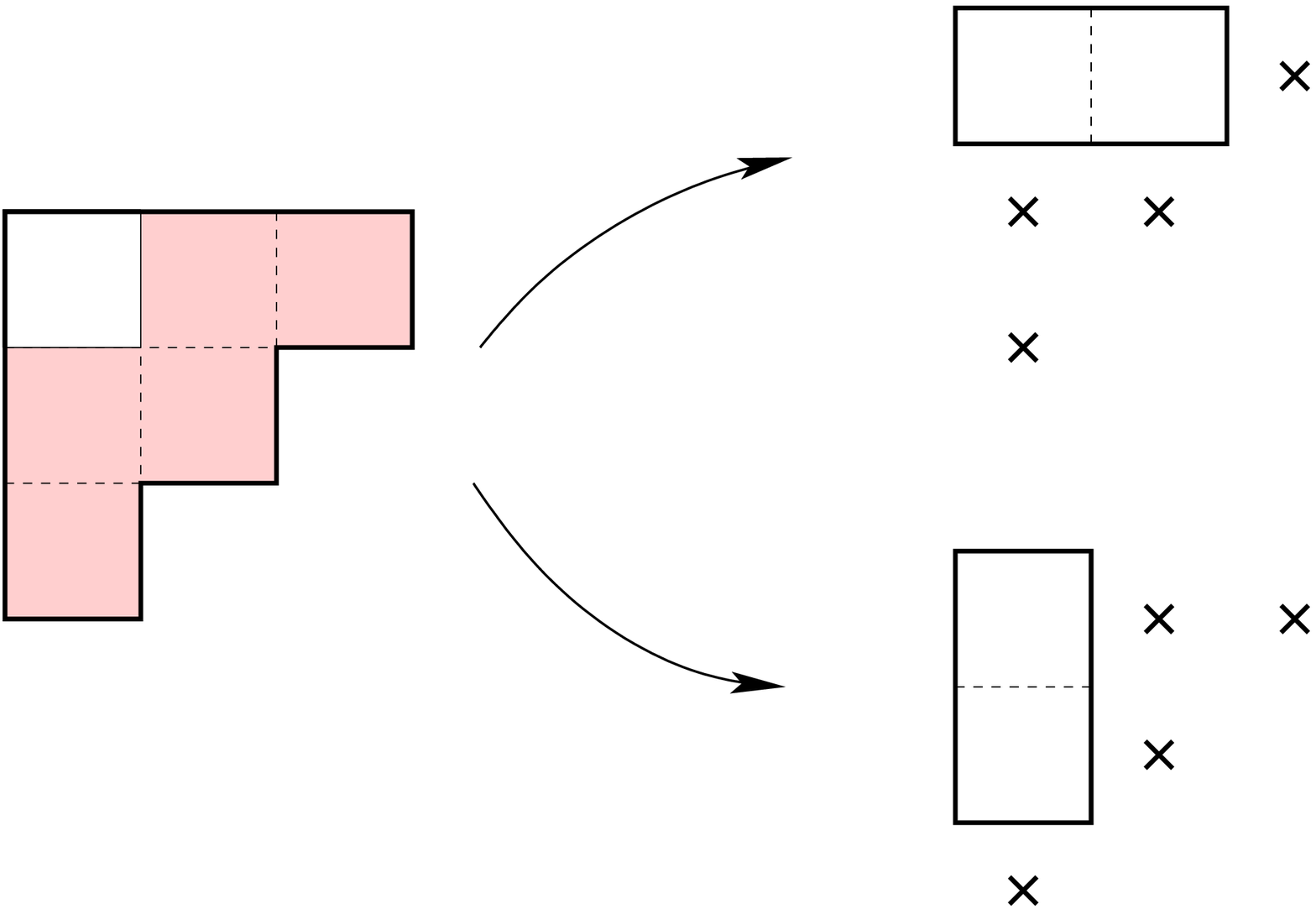}}}
\end{exa}

\begin{proof}[Proof of Theorem \ref{qpieria}] 
By applying the vanishing Corollary \ref{cor2:typea}, we see that it will
suffice to check that the line numbers (that is, the Gromov-Witten invariants
for $d=1$) agree with the second sum in (\ref{pr}). We will sketch the steps
in this argument, and leave the omitted details as an exercise for the reader.

Let $\s_{\ov{\la}}=[X_{\ov{\la}}]$ denote the cohomology class in 
$H^*(G(m+1,E),\Z)$ associated to $\s_{\la}$ for $d=1$, and define
$\s_{\ov{\mu}}$ and $\s_{\ov{p}}$ in a similar way. One then
uses the classical Pieri rule to show that the prescription for the 
line numbers in $\s_{\la}\cdot \s_p$ given in (\ref{pr}) is equivalent
to the identity
\begin{equation}
\label{prequ}
\langle \s_{\la},\s_{\mu}, \s_{p}\rangle_1=\langle
\s_{\ov{\la}}, \s_{\ov{\mu}}, \s_{\ov{p}} \rangle_0,
\end{equation}
where the right hand side of (\ref{prequ}) is a classical intersection 
number on $G(m+1, E)$. 

To prove (\ref{prequ}), observe that the right hand side is given by the
classical Pieri rule on $G(m+1,N)$, and so equals $0$ or $1$. If $\langle
\s_{\ov{\la}}, \s_{\ov{\mu}}, \s_{\ov{p}} \rangle_0=0$, then Corollary
\ref{cor1} shows that $\langle \s_{\la},\s_{\mu}, \s_{p}\rangle_1=0$ as
well. 

Next, assume that $\langle
\s_{\ov{\la}}, \s_{\ov{\mu}}, \s_{\ov{p}} \rangle_0=1$, so that there
is a unique $(m+1)$-dimensional subspace $B$ in the intersection 
$X_{\ov{\la}}(F_\bull)\cap X_{\ov{\mu}}(G_\bull)\cap X_{\ov{p}}(H_\bull)$,
for generally chosen reference flags. Note that the construction of $B$ 
ensures that it lies in the intersection of the corresponding three
Schubert {\em cells} in $G(m+1,E)$, where the defining inequalities 
in (\ref{defequ}) are all equalities. It follows that the two subspaces
\[
V_m=B\cap F_{N-\la_m} \ \ \ \mathrm{and} \ \ \ 
V_m'=B\cap G_{N-\mu_m}
\]
each have dimension $m$, and in fact $V_m\in X_{\la}(F_\bull)$ and 
$V'_m\in X_{\mu}(G_\bull)$. Since 
\[
|\la|+|\mu|=mn+N-p>\dim G(m,N),
\]
we see that $X_\la(F_{\bull})\cap X_\mu(G_\bull)=\emptyset$, and hence
$V_m\neq V'_m$. As $V_m$ and $V_m'$ are both codimension one subspaces of $B$,
this proves that $A=V_m\cap V'_m$ has dimension $m-1$. We deduce that the only
line (corresponding to the required map $f:\bP^1\ra X$ of degree one) meeting
the three Schubert varieties $X_\la(F_{\bull})$, $X_\mu(G_\bull)$, and
$X_p(H_\bull)$ is the locus $\{V\in X\ |\ A\subset V \subset B\}$.
\end{proof}

We conclude with Siebert and Tian's 
presentation of $QH^*(G(m,N))$ in terms of generators and relations.

\begin{thm}[Ring presentation, \cite{ST}] 
\label{qpresa}
The ring $QH^*(X)$ is presented as a quotient of the
polynomial ring $\Z[\s_1,\ldots,\s_n,q]$ by the relations
\[
D_{m+1}=\cdots=D_{N-1}=0  \ \ \ \mathrm{and}  \ \ \ 
D_N+(-1)^nq=0,
\]
where $D_k=\det(\s_{1+j-i})_{1\lequ i,j\lequ k}$ for each $k$.
\end{thm}
\begin{proof} We will justify why the above relations hold in 
$QH^*(X)$, and then sketch the rest of the argument. Since the
degree of $q$ is $N$, the relations $D_k=0$ for $k<N$, which 
hold in $H^*(X,\Z)$, remain true in $QH^*(X)$. For the last relation
we use the formal identity of Schur determinants
\[
D_N-\s_1D_{N-1}+\s_2D_{N-2}-\cdots +(-1)^n\s_n D_m = 0
\]
to deduce that $D_N= (-1)^n\s_nD_m=(-1)^n\s_n\s_{(1^m)}$. Therefore
it will suffice to show that $\s_n\s_{(1^m)}=q$; but this is a 
consequence of Theorem \ref{qpieria}. 

With a bit more work, one
can show that the quotient ring in the theorem is in fact isomorphic
to $QH^*(X)$ (see e.g.\ \cite{Buch}). Alternatively, one may use an
algebraic result of Siebert and Tian \cite{ST}. This
states that for a homogeneous
space $X$, given a presentation
\[
H^*(X,\Z)=\Z[u_1,\ldots,u_r]/(f_1,\ldots,f_t)
\]
of $H^*(X,\Z)$ in terms of homogeneous generators and relations, if
$f_1',\ldots,f_t'$ are homogeneous elements in $\Z[u_1,\ldots,u_r,q]$ such
that $f'_i(u_1,\ldots,u_r,0)=f_i(u_1,\ldots,u_r)$ in $\Z[u_1,\ldots,u_r,q]$
and $f'_i(u_1,\ldots,u_r,q)=0$ in $QH^*(X)$, then the canonical map
\[
\Z[u_1,\ldots,u_r,q]/(f'_1,\ldots,f'_t)\ra QH^*(X)
\]
is an isomorphism. For a proof of this, see \cite[Prop.\ 11]{FP}.
\end{proof}

\medskip
\noin {\bf Remark.} The proofs of Theorems \ref{qgiamba}, \ref{qpieria}, and
\ref{qpresa} do not require the full force of our main Theorem
\ref{thm:typea}.  Indeed, the notion of the kernel and span of a map to $X$
together with Lemma \ref{buchlemma} suffice to obtain the simple proofs
presented here. For instance, to prove Corollary \ref{cor1} one can check
directly that the span of a rational map which contributes to the
Gromov-Witten invariant $\langle \s_{\la}, \s_{\m}, \s_{\nu} \rangle_d$ must
lie in the intersection $X_{\ov{\la}}(F_\bull)\cap X_{\ov{\mu}}(G_\bull)\cap
X_{\ov{\nu}}(H_\bull)$ in $G(m+d,E)$.  
This was the original approach in \cite{Buch}.

The proof of Theorem \ref{qgiamba} used the surprising fact that
in the expansion of the Schur determinant in the quantum Giambelli
formula, each individual monomial is purely classical, that is, has no
$q$ correction terms. This was also observed and generalized to 
partial flag varieties by Ciocan-Fontanine \cite[Thm.\ 3.14]{C-F}.

\newpage

\section{Lecture Four}

\subsection{Schur polynomials} 
\label{sp}
People have known for a long time about the relation between the product of
Schubert classes in the cohomology ring of $G(m,N)$ and the multiplication of
{\em Schur polynomials}, which are the characters of irreducible polynomial
representations of $GL_n$.  Recall that if $Q$ denotes the universal (or
tautological) quotient bundle of rank $n$ over $X$, then the special Schubert
class $\s_i$ is just the $i$th Chern class $c_i(Q)$. If the variables
$x_1,\ldots,x_n$ are the Chern roots of $Q$, then the Giambelli formula
implies that for any partition $\lambda$,
\begin{equation}
\label{basedef}
\s_\lambda =
\det(c_{\la_i+j-i}(Q))= \det(e_{\la_i+j-i}(x_1,\ldots,x_n))= 
s_{\la'}(x_1,\ldots,x_n),
\end{equation}
where $\la'$ is the conjugate partition to $\la$ (whose diagram is the
transpose of the diagram of $\lambda$), and $s_{\la'}(x_1,\ldots,x_n)$ is
a Schur $S$-polynomial in the variables $x_1,\ldots,x_n$. The Schur
polynomials $s_{\lambda}(x_1,\ldots,x_n)$ for $\lambda$ of length at most $n$
form a $\Z$-basis for the ring $\Lambda_n=\Z[x_1,\ldots,x_n]^{S_n}$ of
symmetric polynomials in $n$ variables. It follows that the structure
constants $N_{\la\mu}^{\nu}$ for Schur polynomials
\[
s_{\la}s_{\mu}=\sum_{\nu}N_{\la\mu}^{\nu}s_{\nu}
\]
agree with the Schubert structure constants $c_{\la\mu}^{\nu}$ in 
$H^*(X,\Z)$. 

Our main theorem implies that the Gromov-Witten invariants $\langle \s_{\la},
\s_{\m}, \s_{\n} \rangle_d$ are structure constants in the product of the two
{\em Schubert polynomials} indexed by the permutations for the modified
Schubert varieties $X_{\la}^{(d)}$ and $X_{\mu}^{(d)}$. Postnikov \cite{P} has
shown how one may obtain the same numbers as the coefficients when certain
`toric Schur polynomials' are expanded in the basis of the regular Schur
polynomials.

For the rest of these lectures, we will present the analogue of the
theory developed thus far in the other classical Lie types. To save
time, there will be very little discussion of proofs, but only an
exposition of the main results.  The arguments are often analogous to
the ones in type $A$, but there are also significant differences. For
instance, we shall see that in the case of maximal isotropic
Grassmannians, the equation directly analogous to (\ref{basedef})
defines a family of `$\wt{Q}$-polynomials'. The latter polynomials
have the property that the structure constants in their product
expansions contain both the classical and quantum invariants for these
varieties.

\subsection{The Lagrangian Grassmannian $LG(n,2n)$}
\label{lg}
We begin with the symplectic case and work with the Lagrangian Grassmannian
$LG=LG(n,2n)$ parametrizing Lagrangian subspaces of $E=\C^{2n}$ equipped with
a symplectic form $\langle \, \ ,\ \rangle$.  Recall that a subspace $V$ of
$E$ is {\em isotropic} if the restriction of the form to $V$ vanishes. The
maximal possible dimension of an isotropic subspace is $n$, and in this case
$V$ is called a {\em Lagrangian} subspace. The variety $LG$ is the projective
complex manifold of dimension $n(n+1)/2$ which parametrizes Lagrangian
subspaces in $E$.

The Schubert varieties $X_{\la}(F_\bull)$ in $LG(n,2n)$ now depend on a strict
partition $\lambda=(\la_1>\la_2>\cdots >\la_{\ell}>0)$ with $\la_1\lequ n$;
we let $\D_n$ denote the parameter space of all such $\la$ (a partition
is {\em strict} if all its parts are distinct). We also require a
{\em complete isotropic flag} of subspaces of $E$:
\[
0=F_0\, \subset \, F_1 \, \subset
\cdots \subset F_n\, \subset E
\]
where $\dim(F_i)=i$ for each $i$, and $F_n$ is Lagrangian. 
The codimension $|\la|$
Schubert variety $X_{\la}(F_{\bull})\subset LG$ is defined as the locus of 
$V \in LG$ such that 
\begin{equation}
\label{lgdef}
\dim(V\cap F_{n+1-\la_i})\gequ i,\ \ \mathrm{for} \ \ i=1,\ldots,\ell(\la).
\end{equation}
Let $\s_{\la}$ be the class of $X_{\la}(F_{\bull})$ in the cohomology group
$H^{2|\la|}(LG,\Z)$. We then have a similar list of classical facts, analogous
to those for the type $A$ Grassmannian. However,  
these results were obtained much
more recently than the theorems of Pieri and Giambelli.

\medskip
\noin
{\bf 1)} We have $\dis H^*(LG,\Z)\cong \bigoplus_{\la\in \D_n}\Z\,\s_{\la}$,
that is, the cohomology group of $LG$ is free abelian with basis given 
by the Schubert classes $\s_{\la}$.

\medskip
\noin
{\bf 2)} There is an equation $\s_{\la}\s_{\mu}
=\sum_{\nu}e_{\la\mu}^{\nu}\s_{\nu}$ in $H^*(LG,\Z)$, with
\begin{equation}
\label{agree}
e_{\la\mu}^{\nu}=\int_{LG} \s_{\la}\s_{\mu}\s_{\nu^{\vee}} =  
\# X_{\la}(F_\bull)\cap X_{\mu}(G_\bull) \cap X_{\nu^{\vee}}(H_\bull),
\end{equation}
for general complete isotropic flags 
$F_\bull$, $G_\bull$ and $H_\bull$ in $E$. Here the `dual' partition
$\nu^{\vee}$ is again defined so that $\int_{LG}\s_{\la}\s_{\mu}=
\delta_{\la^{\vee}\mu}$, and it has the property that the parts of 
$\nu^{\vee}$ are the complement of the parts of $\nu$ in the set
$\{1,\ldots,n\}$. For example, the partitions $(4,2,1)$ and $(5,3)$ 
form a dual pair in $\D_5$. 

Stembridge \cite{Ste} has given a combinatorial rule similar to the classical
Littlewood-Richardson rule, which expresses the structure constants
$e_{\la\mu}^{\nu}$ in terms of certain sets of shifted Young tableaux. It
would be interesting to find an analogue of the `puzzle rule' of Theorem
\ref{thm:ktw} that works in this setting.

\medskip
\noin
{\bf 3)} The classes $\s_1,\ldots,\s_n$ are called {\em special Schubert
classes}, and again we have $H^2(LG,\Z)=\Z\,\s_1$. If
\[
0\ra S\ra E_X \ra Q \ra 0
\]
denotes the tautological short exact sequence of vector bundles over $LG$,
then we can use the symplectic form on $E$ to identify $Q$ with the dual of
the vector bundle $S$, and we have $\s_i=c_i(S^*)$, for $0\lequ i \lequ n$.

Let us say that two boxes in a (skew) diagram $\alpha$ are connected if they
share a vertex or an edge; this defines the connected components of $\alpha$.
We now have the following Pieri rule for $LG$, due to Hiller and Boe.

\begin{thm}[Pieri rule for $LG$, \cite{HB}] 
For any $\lambda\in\D_n$ and $p\gequ 0$ we have
\begin{equation}
\label{pieric}
\s_{\lambda} \, \s_p=\sum_{\m} 2^{N(\lambda,\m)}\s_{\m}
\end{equation}
in $H^*(LG, \Z)$, where the sum is over all strict partitions $\mu$
obtained from $\lambda$ by adding $p$ boxes, with no two in the same
column, and $N(\la,\mu)$ is the number of connected components of 
$\mu/\la$ which do not meet the first column.
\end{thm}

\medskip
\noin
{\bf 4)} The Pieri rule (\ref{pieric}) agrees with the analogous product
of Schur $Q$-functions. This was used by Pragacz to obtain a Giambelli 
formula for $LG$, which expresses each Schubert class as a polynomial
in the special Schubert classes.

\begin{thm}[Giambelli formula for $LG$, \cite{P}]
For $i>j>0$, we have 
\[
\s_{i,j}=
\s_i\s_j+2\sum_{k=1}^{n-i}(-1)^k\s_{i+k}\s_{j-k}, 
\]
while for $\la$ of length greater than two,
\begin{equation}
\label{lggiam2}
\s_{\lambda}=\text{\em Pfaffian}[\s_{\lambda_i,\lambda_j}]_{1\lequ i<j\lequ r},
\end{equation}
where $r$ is the smallest even integer such that $r \gequ \ell(\la)$.
\end{thm}

For those who are not so familiar with Pfaffians, we recall that they are 
analogous to (and in fact, square roots of) determinants; see e.g.\ 
\cite[Appendix D]{FPr} for more information. The Pfaffian 
formula (\ref{lggiam2}) is equivalent to the Laplace-type
expansion for Pfaffians
\[
\s_{\la}=\sum_{j=1}^{r-1}(-1)^{j-1}\s_{\la_j,\la_r}
\s_{\la\smallsetminus\{\la_j,\la_r\}}.
\]

\medskip
\noin
{\bf 5)} The ring $H^*(LG,\Z)$ is presented as a quotient of the polynomial
ring $\Z[c(S^*)]=\Z[\s_1,\ldots,\s_n]$ modulo the relations coming from 
the Whitney sum formula
\begin{equation}
\label{pequ}
c_t(S)c_t(S^*)= (1-\s_1t+\s_2t^2-\cdots)(1+\s_1t+\s_2t^2+\cdots)=1.
\end{equation}
By equating the coefficients of like powers of $t$ in (\ref{pequ}), we see that the 
relations are given by
\[
\s_i^2+2\sum_{k=1}^{n-i}(-1)^k\s_{i+k}\s_{i-k}=0
\]
for $1\lequ i\lequ n$, 
where it is understood that $\s_0=1$ and $\s_j=0$ for $j<0$.
In terms of the Chern roots $x_1,\ldots,x_n$ of $S^*$, the equations
(\ref{pequ}) may be written as
\[
\prod_i (1-x_it)\prod_i(1+x_it) = \prod_i(1-x_i^2t^2)=1.
\]
We thus see that $H^*(LG,\Z)$ is isomorphic to the ring 
$\Lambda_n=\Z[x_1,\ldots,x_n]^{S_n}$ modulo the relations 
$e_i(x_1^2,\ldots,x_n^2)=0$, for $1\lequ i\lequ n$.

\subsection{$\wt{Q}$-polynomials} We turn now to the 
analogue of Schur's $S$-polynomials in type $C$, as suggested by
the discussion in \S \ref{sp}. These are a family of polynomials
symmetric in 
the variables $X=(x_1,\ldots,x_n)$, which are modelled on Schur's
$Q$-polynomials. They were defined by Pragacz and Ratajski \cite{PR}
in the course of their work on degeneracy loci.

For strict partitions $\la\in \D_n$, the polynomials 
$\wt{Q}_{\la}(X)$ are obtained by writing
\[
\s_{\la}=\wt{Q}_{\la}(S^*)=\wt{Q}_{\la}(x_1,\ldots,x_n)
\]
as a polynomial in the Chern roots of $S^*$, as we did in (\ref{basedef}).
So $\wt{Q}_i(X)= e_i(X)$ for $0\lequ i\lequ n$, 
\[
\wt{Q}_{i,j}(X)=\wt{Q}_i(X)\wt{Q}_j(X)+2\sum_{k=1}^{n-i}(-1)^k
\wt{Q}_{i+k}(X)\wt{Q}_{j-k}(X),
\]
for $i>j>0$, and for $\ell(\la)\gequ 3$, 
\[
\wt{Q}_{\la}(X)=
\text{Pfaffian}[\wt{Q}_{\lambda_i,\lambda_j}(X)]_{1\lequ i<j\lequ r}.
\]
Pragacz and Ratajski noticed that this definition also makes sense for {\em
non-strict} partitions $\lambda$. For $\la_1> n$, one checks easily that
$\wt{Q}_{\la}(X)=0$. Let $\E_n$ denote the parameter space of all partitions
$\lambda$ with $\lambda_1\lequ n$. We then obtain polynomials
$\wt{Q}_{\la}(X)$ for $\la\in \E_n$ with the following properties:

\medskip
\noin
a) The set $\{\wt{Q}_{\la}(X)\ |\ \la\in \E_n\}$ is a free $\Z$-basis for
$\Lambda_n$. 

\medskip
\noin
b) $\wt{Q}_{i,i}(X)=e_i(x_1^2,\ldots,x_n^2)$, for $1\lequ i\lequ n$.

\medskip
\noin c) (Factorization Property) If $\la=(\la_1,\ldots,\la_{\ell})$ and
$\la^+$ is defined by $\la^+=\la\cup
(i,i)=(\la_1,\ldots,i,i,\ldots,\la_{\ell})$, then
\[
\wt{Q}_{\la^+}(X)=\wt{Q}_{\la}(X)\cdot\wt{Q}_{i,i}(X).
\]

\medskip
\noin
d) For strict $\la$, the $\wt{Q}_{\la}(X)$ enjoy the same Pieri rule 
as in (\ref{pieric})
\begin{equation}
\label{pieriwtq}
\wt{Q}_{\lambda}(X)\cdot\wt{Q}_p(X)=\sum_{\m} 
2^{N(\lambda,\m)}\wt{Q}_{\m}(X),
\end{equation}
only this time the sum in (\ref{pieriwtq}) is over all partitions $\mu\in
\E_n$ (possibly not strict) obtained from $\la$ by adding $p$ boxes, with no
two in the same column.  In particular, it follows that
\[
\wt{Q}_n(X)\cdot\wt{Q}_{\la}(X)= \wt{Q}_{(n,\la)}(X)
\]
for all $\la\in \E_n$.

\medskip
\noin
e) There are structure constants $e_{\la\mu}^{\nu}$ such that
\[
\wt{Q}_{\la}(X)\cdot\wt{Q}_{\mu}(X)=\sum_{\nu}
e_{\la\mu}^{\nu}\,\wt{Q}_{\nu}(X),
\]
defined for $\la,\mu,\nu\in\E_n$ with $|\nu|=|\la|+|\mu|$. These agree
with the integers in (\ref{agree}) if $\la$, $\mu$, and $\nu$ are 
strict. In general, however, these integers can be negative, for example
\[
e_{(3,2,1),(3,2,1)}^{(4,4,2,2)}=-4.
\]
In the next lecture, 
we will see that some of the constants $e_{\la\mu}^{\nu}$ for non-strict $\nu$
must be positive, as they
are equal to three-point Gromov-Witten invariants, up to a power of $2$.

\medskip
Finally, observe that the above properties allow us to present the cohomology
ring of $LG(n,2n)$ as the quotient of the ring $\Lambda_n=\Z[X]^{S_n}$ of
$\wt{Q}$-polynomials in $X$ modulo the relations $\wt{Q}_{i,i}(X)=0$, for
$1\lequ i\lequ n$.

\newpage

\section{Lecture Five}

\subsection{Gromov-Witten invariants on $LG$} 
As in the first lecture, by 
a rational map to $LG$ we mean a morphism $f\colon \bP^1\to LG$, and its
degree is the degree of $f_*[\bP^1]\cdot\s_1$. The Gromov-Witten invariant
$\langle \s_\lambda, \s_\mu, \s_\nu\rangle_d$ is defined for
$|\lambda|+|\mu|+|\nu|=n(n+1)/2+d(n+1)$ and counts the number of
rational maps $f\colon\bP^1\ra LG(n,2n)$ of degree $d$ such that 
$f(0)\in X_\lambda(F_\bull)$, $f(1)\in X_\mu(G_\bull)$,
and $f(\infty)\in X_\nu(H_\bull)$, for given flags $F_\bull$,
$G_\bull$, and $H_\bull$ in general position.

We also define the kernel of a map $f:\bP^1\ra LG$ as the intersection of the
subspaces $f(p)$ for all $p\in \bP^1$. In the symplectic case it happens that
the span of $f$ is the orthogonal complement of the kernel of $f$, and hence
is not necessary. Therefore the relevant parameter space of kernels that
replaces the two-step flag variety is the isotropic Grassmannian $IG(n-d,2n)$,
whose points correspond to isotropic subspaces of $E$ of dimension $n-d$.

If $d\gequ 0$ is an integer, $\lambda$, $\mu$, $\nu\in\D_n$ are
such that $|\lambda|+|\mu|+|\nu|=n(n+1)/2+d(n+1)$, 
and $F_\bull$, $G_\bull$, and $H_\bull$ are complete isotropic flags
of $E=\C^{2n}$ in general position, then similar arguments to the ones
discussed earlier show that 
the map $f\mapsto\Ker(f)$ gives a bijection of the set of
rational maps $f\colon \bP^1 \ra LG$ of degree $d$ satisfying
$f(0)\in X_\lambda(F_\bull)$, $f(1)\in X_\mu(G_\bull)$, and
$f(\infty)\in X_\nu(H_\bull)$, with the set of points 
in the intersection
$X^{(d)}_\lambda(F_\bull) \cap X^{(d)}_\mu(G_\bull) \cap
X^{(d)}_\nu(H_\bull)$ in $Y_d=IG(n-d,2n)$. We therefore get

\begin{cor}[\cite{BKT1}]
\label{cor:typec}
Let $d\gequ 0$ and
$\la$, $\m$, $\n\in\D_n$ be chosen as above. Then 
\[
\langle \s_\la, \s_\m, \s_\n \rangle_d=
\int_{IG(n-d,2n)} [X^{(d)}_{\la}]\cdot [X^{(d)}_{\m}]\cdot [X^{(d)}_{\nu}].
\] 
\end{cor}

The {\em line numbers} $\langle \sigma_{\la},\sigma_\m,\sigma_\n\rangle_1$
satisfy an additional relation, which is an extra ingredient needed
to complete the analysis for $LG(n,2n)$.

\begin{prop}[\cite{KTlg}]
\label{linenumbers}
For $\la$, $\m$, $\n\in \D_n$ we have 
\[
\langle \sigma_{\la},\sigma_\m,\sigma_\n\rangle_1=
\frac{1}{2}\int_{LG(n+1,2n+2)}[X^+_\la]\cdot[X^+_\m]\cdot[X^+_\n],
\]
where $X_{\la}^+$, $X^+_\m$, $X^+_\n$ 
 denote Schubert varieties in $LG(n+1,2n+2)$.
\end{prop}

The proof of Proposition \ref{linenumbers} in \cite{KTlg} proceeds 
geometrically, by using a correspondence between lines on
$LG(n,2n)$ (which are parametrized by points of $IG(n-1,2n)$) and
points on $LG(n+1,2n+2)$. 

\subsection{Quantum cohomology of $LG(n,2n)$} 
\label{lgapp}
The quantum cohomology ring of $LG$ is a $\Z[q]$-algebra isomorphic to
$H^*(LG,\Z)\otimes_{\Z}\Z[q]$ as a module over $\Z[q]$, but here $q$ is a
formal variable of degree $n+1$. The product in $QH^*(LG)$ is defined by the
same equation (\ref{qprod}) as before, but as  $\deg(q)=n+1<2n$, we expect
different behavior than what we have seen for $G(m,N)$. The
previous results allow one to prove the following theorem (the
original proofs in \cite{KTlg} were more involved).

\begin{thm}[Ring presentation and quantum Giambelli, \cite{KTlg}]
\label{lggiambthm} 
The ring $QH^*(LG)$ is presented as a quotient of the polynomial
ring $\Z[\s_1,\ldots,\s_n,q]$ by the relations 
\[
\s_i^2+2\sum_{k=1}^{n-i}(-1)^k\s_{i+k}\s_{i-k}=(-1)^{n-i}\s_{2i-n-1}\,q
\]
for $1\lequ i\lequ n$. The Schubert class $\s_{\lambda}$
in this presentation is given by the Giambelli formulas
\[
\s_{i,j}=
\s_i\s_j+2\sum_{k=1}^{n-i}(-1)^k\s_{i+k}\s_{j-k}+(-1)^{n+1-i}\s_{i+j-n-1}\,q
\]
for $i>j>0$, and for $\ell(\la)\gequ 3$,
\begin{equation}
\label{lgqgiam2}
\s_{\lambda}=\text{\em Pfaffian}[\s_{\lambda_i,\lambda_j}]_{1\lequ i<j\lequ r}.
\end{equation}
\end{thm}

The key observation here is that the quantum Giambelli formulas for $LG(n,2n)$
coincide with the classical Giambelli formulas for $LG(n+1,2n+2)$, when the
class $2\s_{n+1}$ is identified with $q$. Note that this does not imply that
the cohomology ring of $LG(n+1,2n+2)$ is isomorphic to $QH^*(LG(n,2n))$,
because the relation $\s_{n+1}^2=0$, which holds in the former ring, does not
hold in the latter (as $q^2\neq 0$).

Using the $\wt{Q}$-polynomials, we can write the presentation of $QH^*(LG)$ as
follows. Let $X^+=(x_1,\ldots,x_{n+1})$ and let $\wt{\Lambda}_{n+1}$ be the
subring of $\Lambda_{n+1}$ generated by the polynomials $\wt{Q}_i(X^+)$ for
$1\lequ i\lequ n$ together with $2\,\wt{Q}_{n+1}(X^+)$. Then the map
$\wt{\Lambda}_{n+1}\ra QH^*(LG)$ which sends $\wt{Q}_{\la}(X^+)$ to $\s_{\la}$
for $\la\in \D_n$ and $2\,\wt{Q}_{n+1}(X^+)$ to $q$ extends to a surjective
ring homomorphism, whose kernel is generated by the relations
$\wt{Q}_{i,i}(X^+)=0$ for $1\lequ i\lequ n$.

Theorem \ref{lggiambthm} therefore implies that the algebra in
$QH^*(LG)$ is controlled by the multiplication of $\wt{Q}$-polynomials.
In particular, the quantum Pieri rule for $LG$ is a specialization 
of the Pieri rule for $\wt{Q}$-polynomials.

\begin{thm}[Quantum Pieri rule for $LG$, \cite{KTlg}] 
\label{qupieri}
For any $\lambda\in\D_n$ and $p\gequ 0$ we have
\begin{equation}
\label{quantumpieri}
\s_{\lambda}\cdot\s_p=\sum_{\m} 2^{N(\lambda,\m)}\s_{\m}+\sum_{\n}
2^{N'(\nu,\lambda)} \s_{\nu} \, q
\end{equation}
in $QH^*(LG(n,2n))$, where the first sum is classical, as in (\ref{pieric}),
while the second is over all strict $\n$ obtained from $\lambda$ by
subtracting $n+1-p$ boxes, no two in the same column, and $N'(\nu,\lambda)$ is
one less than the number of connected components of $\la/\nu$.
\end{thm}
 
The informed reader will notice that the exponents $N'(\nu,\lambda)$ in the
multiplicities of the quantum correction terms in (\ref{quantumpieri}) are of
the kind encountered in the classical Pieri rule for {\em orthogonal}
Grassmannians. This was the first indication of a more general phenomenon,
which we will discuss at the end of this lecture.

For arbitrary products in $QH^*(LG)$, we have

\begin{cor}
\label{ancor} 
In the relation 
\[
\s_{\la}\cdot\s_{\mu}=\sum_{\substack{d\gequ 0 \\ |\nu|=|\la|+|\mu|-
d(n+1)}} \langle \s_{\la}, \s_{\mu},\s_{\nu^{\vee}}\rangle_d\,\s_{\nu}\, q^d, 
\]
the quantum structure constant 
$\langle \s_{\la}, \s_{\mu},\s_{\nu^{\vee}}\rangle_d$ is equal to 
$2^{-d}e_{\la,\mu}^{((n+1)^d,\nu)}$.
\end{cor}

Corollary \ref{ancor} follows immediately from Theorem \ref{lggiambthm}
together with the identity
\[
\wt{Q}_{((n+1)^d,\nu)}(X^+)=\wt{Q}_{n+1}(X^+)^d\cdot \wt{Q}_{\lambda}(X^+)
\]
of $\wt{Q}$-polynomials. We deduce that the $\wt{Q}$-polynomial structure
constants of the form $e_{\la,\mu}^{(n^d,\nu)}$ (for $\la,\mu,\nu\in\D_{n-1}$)
are nonnegative integers, divisible by $2^d$. A combinatorial rule for these
numbers would give a quantum Littlewood-Richardson rule for $LG$.

The proofs of Theorems \ref{lggiambthm} and \ref{qupieri}, as compared to
those for the type $A$ Grassmannian $G(m,N)$, are complicated by two
facts. First, a different argument is needed to establish the quantum Pieri
rule, which is related by Proposition \ref{linenumbers} to the classical Pieri
rule on $LG(n+1,2n+2)$. Second, in the quantum Giambelli Pfaffian expansion
$\s_{\lambda}=\text{Pfaffian} [\s_{\lambda_i,\lambda_j}]_{1\lequ i<j\lequ r}$,
there are terms which do involve $q$-corrections, and these extra $q$-terms
cancel each other out in the end. Thus more combinatorial work is required to
prove that the Pfaffian formula (\ref{lgqgiam2}) holds in $QH^*(LG)$.

\subsection{The orthogonal Grassmannian $OG(n+1,2n+2)$}
\label{og}
We now turn to the analogue of the above theory in the orthogonal Lie
types. We will work with the even orthogonal Grassmannian
$OG=OG(n+1,2n+2)=SO_{2n+2}/P_{n+1}$. This variety parametrizes (one component
of) the locus of maximal isotropic subspaces of a $(2n+2)$-dimensional vector
space $E$, equipped with a nondegenerate symmetric form. Note that there are
two families of such subspaces; by convention, given a fixed isotropic flag
$F_\bull$ in $E$, we consider only those isotropic $V$ in $E$ such that $V\cap
F_{n+1}$ has even codimension in $F_{n+1}$. We remark that $OG$ is isomorphic
to the odd orthogonal Grassmannian $OG(n,2n+1)=SO_{2n+1}/P_n$, hence our
analysis (for the maximal isotropic case) will include both the Lie types $B$
and $D$.

The Schubert varieties $X_{\lambda}(F_{\bull})$ in $OG$ are again parametrized
by partitions $\lambda\in\D_n$ and are defined by the same equations
(\ref{lgdef}) as before, with respect to a complete isotropic flag $F_{\bull}$
in $E$.  Let $\tau_{\lambda}$ be the cohomology class of $X_{\la}(F_{\bull})$;
the set $\{\tau_{\la}\, |\, \la\in\D_n\}$ is a $\Z$-basis of $H^*(OG,\Z)$.
Now much of the theory for $OG$ is similar to that for $LG(n,2n)$. To save
time, we will pass immediately to the results about the quantum cohomology
ring of $OG$. We again have an isomorphism of $\Z[q]$-modules $QH^*(OG)\cong
H^*(OG,\Z)\otimes\Z[q]$, but this time the variable $q$ has degree $2n$.

Another difference between the symplectic and orthogonal case is that the
natural embedding of $OG(n+1,2n+2)$ into the type $A$ Grassmannian
$G(n+1,2n+2)$ is degree doubling. This means that for every degree $d$ map
$f:\bP^1\ra OG$, the pullback of the tautological quotient bundle over $OG$
has degree $2d$. It follows that the relevant parameter space of kernels of
the maps counted by a Gromov-Witten invariant is the sub-maximal isotropic
Grassmannian $OG(n+1-2d,2n+2)$. We pass directly to the corollary of the
corresponding `main theorem':

\begin{cor}[\cite{BKT1}]
\label{cor:typeb}
Let $d\gequ 0$ and
$\la$, $\m$, $\n\in\D_n$ be such that $|\la|+|\mu|+|\nu|=n(n+1)/2+2nd$.
Then
\[
\langle \tau_\la, \tau_\m, \tau_\n \rangle_d=
\int_{OG(n+1-2d,2n+2)} 
[X^{(d)}_{\la}]\cdot [X^{(d)}_{\m}]\cdot [X^{(d)}_{\nu}].
\] 
\end{cor}

This result may be used to obtain analogous structure theorems for $QH^*(OG)$.

\begin{thm}[Ring presentation and quantum Giambelli,\cite{KTorth}]
\label{oggiambthm} 
The ring $QH^*(OG)$ is presented as a quotient of the polynomial
ring $\Z[\tau_1,\ldots,\tau_n,q]$ modulo the relations 
\[
\tau_i^2+2\sum_{k=1}^{i-1}(-1)^k\tau_{i+k}\tau_{i-k}+(-1)^i\tau_{2i}=0
\]
for all $i<n$, together with the quantum relation
\[
\tau_n^2=q.
\]
The Schubert class $\tau_{\lambda}$
in this presentation is given by the Giambelli formulas
\[
\tau_{i,j}=
\tau_i\tau_j+2\sum_{k=1}^{j-1}(-1)^k\tau_{i+k}\tau_{j-k}+(-1)^j\tau_{i+j}
\]
for $i>j>0$, and for $\ell(\la)\gequ 3$,
\[
\tau_{\lambda}=\text{\em Pfaffian}
[\tau_{\lambda_i,\lambda_j}]_{1\lequ i<j\lequ r}.
\]
\end{thm}

It follows from this that the quantum Giambelli formula for $OG$ coincides
with the classical Giambelli formula, and indeed the results in the 
orthogonal case are formally closer to those in type $A$. 

One can use the $\wt{P}$-polynomials, defined by 
$\wt{P}_{\la}=2^{-\ell(\la)}\wt{Q}_{\la}$ for each $\la$, to describe
the multiplicative structure of $QH^*(OG)$. Let $\Lambda'_n$ denote the
$\Z$-algebra generated by the polynomials $\wt{P}_{\la}(X)$, for $\la\in
\E_n$, where $X=(x_1,\ldots,x_n)$. Then the map which sends $\wt{P}_{\la}(X)$
to $\tau_{\la}$ for all $\la\in \D_n$ and $\wt{P}_{n,n}(X)$ to $q$ extends
to a surjective ring homomorphism $\Lambda_n'\ra QH^*(OG)$ with kernel
generated by the relations $\wt{P}_{i,i}(X)=0$, for all $i<n$. 
Define the structure constants $f_{\la\mu}^{\nu}$ by the relation
\[
\wt{P}_{\la}(X)\cdot\wt{P}_{\mu}(X)=\sum_{|\nu|=|\la|+|\mu|}f_{\la\mu}^{\nu}
\, \wt{P}_{\nu}(X).
\]

\begin{cor}
The Gromov-Witten invariant (and quantum structure constant)
$\langle \tau_{\la}, \tau_{\mu},\tau_{\nu^{\vee}}\rangle_d$ is equal to 
$f_{\la,\mu}^{(n^{2d},\nu)}$.
\end{cor}

\begin{thm}[Quantum Pieri rule for $OG$, \cite{KTorth}] 
\label{ogqupieri}
For any $\lambda\in\D_n$ and $p\gequ 0$ we have
\[
\tau_{\lambda}\cdot\tau_p=\sum_{\m} 2^{N'(\lambda,\m)}\tau_{\m}+\
\sum_{\nu}2^{N'(\lambda,\nu)} \tau_{\nu\ssm (n,n)} \, q,
\]
where the first sum is over strict $\mu$ and the second over partitions
$\nu=(n,n,\ov{\nu})$ with $\ov{\nu}$ strict, such that both $\mu$ and
$\nu$ are obtained from $\lambda$ by adding $p$ boxes, with no two in
the same column.
\end{thm}

The above quantum Pieri rule implies that 
\[
\tau_{\la}\cdot\tau_n=
\begin{cases}
\tau_{(n,\la)}&\text{if $\la_1<n$}, 
\\ \tau_{\la\ssm (n)}\,q&\text{if $\la_1=n$}
\end{cases}
\]
in the quantum cohomology ring of $OG(n+1,2n+2)$. We thus see that
multiplication by $\tau_n$ is straightforward; it follows that 
to compute all the Gromov--Witten invariants for
$OG$, it suffices to evaluate the
$\langle \tau_{\la}, \tau_{\m}, \tau_{\n} \rangle_d$ for $\m,\n\in
\D_{n-1}$. Define a map $\wh{{}}:\D_n\ra\D_{n-1}$ by setting
$\wh{\la}=(n-\la_{\ell},\ldots,n-\la_1)$ for
any partition $\la$ of length $\ell$. Notice that \, $\wh{{}}$ \, is 
essentially a type $A$ Poincar\'e duality map (for the type $A$ 
Grassmannian $G(\ell,n+\ell)$).

Partitions in $\D_{n-1}$ also parametrize the Schubert classes $\s_{\la}$
in the cohomology of $LG(n-1,2n-2)$. The two spaces $OG(n+1,2n+2)$ and
$LG(n-1,2n-2)$ have different dimensions and seemingly little in common
besides perhaps the fact the the degree of $q$ in the quantum cohomology of
the former is twice the degree of $q$ in the latter. However, if
${}^\vee:\D_{n-1}\ra\D_{n-1}$ denotes the Poincar\'e duality involution on
$\D_{n-1}$, we have the following result.

\begin{thm}[\cite{KTorth}]
\label{OGLG}
Suppose that $\la\in\D_n$ is a non-zero
partition with $\ell(\la)=2d+e+1$
for some nonnegative integers $d$ and $e$. For any
$\m,\n\in\D_{n-1}$, we have an equality
\begin{equation}
\label{lgogequ}
\langle \tau_{\la}, \tau_{\m}, \tau_{\n} \rangle_d =
\langle \s_{\wh{\la}}, \s_{\m^{\vee}}, \s_{\n^{\vee}} \rangle_e
\end{equation}
of Gromov--Witten
invariants for  $OG(n+1,2n+2)$
and $LG(n-1,2n-2)$, respectively. If $\la$ is zero or
$\ell(\la)<2d+1$, then
$\langle \tau_{\la}, \tau_{\m}, \tau_{\n} \rangle_d =0$.
\end{thm}

We remark that the left hand side of (\ref{lgogequ}) is symmetric
in $\la$, $\mu$, and $\nu$, unlike the right hand side. This reflects 
a $(\Z/2\Z)^3$-symmetry shared by the Gromov-Witten invariants for both
$LG$ and $OG$. In fact, Theorem \ref{OGLG} is essentially equivalent to 
this symmetry. The proof in \cite{KTlg}\cite{KTorth} proceeds by first 
establishing the symmetry by a clever use of the quantum Pieri rule, and 
then using the relation between the structure constants of $\wt{Q}$- and 
$\wt{P}$-polynomials to put everything together. As of this writing, we
lack a purely geometric result that would explain Theorem \ref{OGLG}.

\subsection{Concluding remarks} 
The kernel and span ideas in these notes have been used by 
Buch in \cite{qmonk} and \cite{Bu3} to obtain simple proofs of the 
main structure theorems regarding the quantum
cohomology of any partial flag manifold $SL_N/P$, where $P$ is a 
parabolic subgroup of $SL_N$ (the arguments assume the associativity of
the quantum product). However, in \cite{BKT1} it is shown that 
there is no direct analogue of Theorem \ref{thm:typea} in this
generality, at
least not for the complete flag manifold $SL_N/B$. 

In recent work with Buch and Kresch \cite{BKT2}, we present similar results to
the ones described here (including an analogue of Theorem \ref{thm:typea}) for
any homogeneous space of the form $G/P$, were $G$ is a classical Lie group and
$P$ is a maximal parabolic subgroup of $G$. These manifolds include the
Grassmannians parametrizing non-maximal isotropic subspaces which were
mentioned earlier.

For the reader who is interested in learning more about the classical and
quantum cohomology of homogeneous spaces, we recommend the texts by Fulton
\cite{F} and Manivel \cite{M} and the expository article \cite{FP}. The latter
reference features the general approach to Gromov-Witten theory using
Kontsevich's moduli space of stable maps.

\skipline
\noin
{\bf Acknowledgements.} It is a pleasure 
to thank Piotr Pragacz and Andrzej Weber for their efforts
in organizing the mini-school on `Schubert Varieties' at the Banach
Center in Warsaw. I am also grateful 
to Marc Levine and the Universit\"at Essen for their hospitality and
the Wolfgang Paul program of the Humbolt Foundation for support during 
the period when these notes were written. Finally, I thank my collaborators
Anders Buch and Andrew Kresch for many enlightening discussions, and for 
all the fun we've had working together on this project.

\end{document}